\documentstyle[12pt]{article}
\newcommand{\const}{\mathop{\rm const}\limits}

\newcommand{\mod}{\mathop{\rm mod}\limits}

\newcommand{\supp}{\mathop{\rm supp}\limits}

\newcommand{\vraisup}{\mathop{\rm vraisup}\limits}

\begin{document}

\begin{center}

{\bf HARDY'S OPERATOR AND NORMABILITY OF }\par
\vspace{3mm}

{\bf GENERALIZED LORENTZ -  MARCINKIEWICZ SPACES,  }\par

\vspace{3mm}

{\bf with sharp or weakly sharp constant estimation.} \par

\vspace{4mm}

 $ {\bf E.Ostrovsky^a, \ \ L.Sirota^b } $ \\

\vspace{4mm}

$ ^a $ Corresponding Author. Department of Mathematics and Computer Science, Bar-Ilan University, 84105, Ramat Gan, Israel.\\
\end{center}
E - mail: \ galo@list.ru \  eugostrovsky@list.ru\\
\begin{center}
$ ^b $  Department of Mathematics and computer science. Bar-Ilan University,
84105, Ramat Gan, Israel.\\

E - mail: \ sirota3@bezeqint.net\\

\vspace{3mm}
                    {\sc Abstract.}\\

 \end{center}

 \vspace{3mm}

 We introduce a  Banach rearrangement invariant (tail) quasy - norm by means of Hardy's (Cesaro) average
on the (measurable) functions defined on some measurable space which is a slight generalization
of classical Lorentz - Marcinkiewicz  norm  and  find for it an equivalent norm expression.\par

 \vspace{3mm}

{\it Key words and phrases:} Tail function, rearrangement invariant norm and spaces, weight, fundamental  function,
random variable, Hardy (Cesaro) operator, slowly and regular varying functions, examples, integral operator, factorization,
upper and lower estimates, right and left inverse function, ordinary and Grand Lebesgue  spaces, Lorentz, Marcinkiewicz
norm and spaces, exactness and weak exactness. \par

\vspace{3mm}

 {\it Mathematics Subject Classification (2000):} primary 60G17; \ secondary 60E07; 60G70.\\

\hfill $\Box$ \\
 \bigskip

\section{Notations. Statement of problem.}

\vspace{3mm}

 Let $ (X = \{x\}, \cal{A}, \mu) $ be measurable space with  non-trivial sigma-finite measure $ \mu. $
We will suppose without loss of generality in the case  $ \mu(X) < \infty $  that $ \mu(X) = 1 $ (the probabilistic case)
and denote $  x = \omega, \ {\bf P} = \mu. $\par
 Define as usually for arbitrary measurable function $ f: X \to R $ its distribution function
 (more exactly, tail function)

 $$
 T_f(t) = \mu\{x: |f(x)| \ge t \}, \  t \ge 0,
 $$

 $$
 ||f||_p = \left[ \int_X |f(x)|^p \ \mu(dx) \right]^{1/p}, \ p \ge 1; \ L_p = L_{p,\mu} = \{f, ||f||_p < \infty\},
 $$

$$
||f||_{\infty} = \vraisup_x |f(x)| \ (\mod \mu),
$$
and denote by $ f^*(t) = T_f^{-1}(t) $ the left inverse to the tail function $ T_f(t);  $
$$
 f^{**}(t) \stackrel{def}{=} t^{-1} \int_0^t f^*(s) \ ds, \ t > 0.
$$
 It will be presumed  in the case $ X \subset R $  that the measure $  \mu $  is classical Lebesgue measure:
$ \mu(dx) = dx. $ \par

  We will denote the set of all tail functions  as $ \{ T \}; $ obviously, the set  $ \{ T \}$ contains on
all the functions  which are right continuous, monotonically non-increasing
with values in the set  $ [0, \mu(X)]. $ \par

  Let also $ (V, ||\cdot||V) $ be  arbitrary Banach complete function space  over the functions defined on the
open semi-axis $  R_+ =(0,\infty), $ not necessary to be rearrangement invariant.  This imply in particular that
if $ f_1, f_2 \in V  $ and $ |f_1(t)| \le |f_2(t)|, $ then $ ||f_1||V \le ||f_2||V.  $ \par

\vspace{3mm}
{\bf Definition 1.1.} \par
\vspace{3mm}
{\it We  introduce the following quasy - norms} $ ||| f |||^*_Y= ||| f |||^*_{Y,V} $ {\it and}
$ ||| f |||_Y= ||| f |||_{Y,V} $ {\it for the measurable  functions} $ f: X \to R $ {\it as follows:}

$$
||| f |||^*_Y= ||| f |||^*_{Y,V} \stackrel{def}{=} ||f^*||V,\eqno(1.1)
$$

$$
||| f |||_Y= ||| f |||_{Y,V} \stackrel{def}{=} ||f^{**}||V, \eqno(1.2)
$$
{\it and correspondingly the following spaces}
$$
(Y, |||\cdot|||Y) = \{f, |||f|||_Y < \infty  \}, \hspace{5mm} (Y_*, |||\cdot|||Y_*) = \{f, |||f|||^*_Y < \infty  \}. \eqno(1.3)
$$
\vspace{3mm}

 These spaces (or similar) are  named by D.E.Edmund and B.Opic
 in \cite{Edmund2} as "Lorentz-Karamata spaces", by B.Opic and L.Pick
 in \cite{Opic1} as "Lorentz-Zygmund spaces", by  E.Pustylnik
in \cite{Pustylnik1} as "Ultrasymmetric  spaces".\par

\vspace{3mm}

{\bf  We will prove in this article that under some simple conditions the space $ (Y, |||\cdot|||Y) $ is (complete)
rearrangement invariant Banach space and that the quasy - norm $ ||| f |||^*_Y $ and the norm $ ||| f |||_Y $ are linear
 equivalent:}

$$
K_1(V)\ ||| f |||^*_Y   \le ||| f |||_Y \le  K_2(V) \ ||| f |||^*_Y, \eqno(1.4)
$$
{\bf where $ K_1(V), \  K_2(V)  $ are finite positive constants (more exactly, function on $ V ) $
depending only on the space $ (V, ||\cdot||V) $ but not on the
function $ f, $ and will find sharp (exact) values or as a minimum weak sharp (i.e. up to multiplicative constant) values of
these functions. }\par
 In the articles  \cite{Ostrovsky111}, \cite{Ostrovsky113} the inequalities of a view (1.4) were applied in the
theory of Probability and further - in Statistics  and in the Monte - Carlo method in order to characterize the
tail behavior of random variables and sums of random variables.\par

\vspace{3mm}

 This problem in less general statement see in \cite{Carro1}, \cite{Carro22}, \cite{Cwikel33},  \cite{Edmund1},
 \cite{Maligranda11}, \cite{Maz'ja1}, \cite{Milman1}, \cite{Kaminska1}, \cite{Martin1}, \cite{Milman2}, \cite{Milman3},
\cite{Ostrovsky111}, \cite{Ostrovsky113}, \cite{Soria1}; see also reference therein. \par

 The another but near statement of problem based on the Cesaro average and Cesaro spaces see in the recent article of
S.V.Astashkin  and L.Maligranda  \cite{Astashkin2}. \par
\vspace{3mm}

{\bf Example 1.1.}  (See \cite{Ostrovsky111}).  \par

 Let  $ w = w(s), s \ge 0 $ be any continuous strictly increasing   numerical function (weight) defined on the set
 $ s \in (0, \infty) $  such that
 $$
 w(s) = 0 \Leftrightarrow s = 0;  \ \lim_{s \to \infty}w(s) = \infty. \eqno(1.5)
 $$

   We impose here on the set of all such a functions $ W = \{ w \} $ the following  restriction: \par

$$
\forall w \in W \ \exists T \in \{ T \}  \Rightarrow w(T(s)) = 1/s.
$$

 Let us introduce the following important functional

 $$
 \gamma(w) = \sup_{t > 0} \left[  \frac{w(t)}{t} \ \int_0^t \frac{du}{w(u)} \right]\eqno(1.6)
 $$
and the following quasi-norms:

$$
||f||_w^* = \sup_{t > 0} [w(t) \ f^*(t)], \eqno(1.7)
$$

$$
||f||_w = \sup_{t > 0} [w(t) \ f^{**}(t)],\eqno(1.8)
$$
 The necessary and sufficient condition for finiteness of the functional $ \gamma(w) $ see, e.g.
in the article \cite{Astashkin1}. \par

\vspace{3mm}
{\bf Remark 1.1.} Note that
\vspace{3mm}
 $$
 ||f||^*_w = \sup_{t > 0} [t w(T_f(t))],
 $$
so that if $ ||f||^*_w \in (0,\infty),  $ then
$$
T_f(t) \le w^{-1}(||f||^*_w/t).
$$

 Therefore the functional $ f \to ||f||^* $ may called "the tail quasinorm". \par

\vspace{3mm}
{\bf Remark 1.2.} As long as
\vspace{3mm}
$$
f^{**}(t)= t^{-1} \sup_{\mu(E) \le t} \int_E |f(x)| \ \mu(dx),\eqno(1.9)
$$
we can rewrite  the expression for $ ||f||_w  $ as follows:

$$
||f||_w = \sup_{t > 0} \left[ (w(t)/t) \cdot \sup_{E: \mu(E) \le t} \int_E  |f(x)| \ \mu(dx)   \right]. \eqno(1.10)
$$

 If the measure $ \mu $ has not atoms, then the expression (1.10) may be rewritten as follows:

 $$
||f||_w = \sup_{ E: 0 < \mu(E) < \infty } \left[ \frac{ w(\mu(E))}{\mu(E)} \cdot \int_E |f(x)| \ \mu(dx) \right]. \eqno(1.11)
 $$

  It follows from  equality (1.10) that $ ||f||_w $ is natural rearrangement invariant norm and the space
$ L_w = \{ f: ||f||_w < \infty  \}  $ is complete Banach functional rearrangement invariant
 space with Fatou property.  The proof is similar to one in the case $  w(t) = t^{1/p}, \ p \ge 1;  $ see
\cite{Bennet1}, chapters 1,2; \cite{Stein1}, chapter 8. \par

\vspace{3mm}

 The norm $  ||f||_w   $ is named Marcinkiewicz's norm, see \cite{Krein1},  chapter 2, section 2. \par

\vspace{3mm}
 It is proved in \cite{Ostrovsky111} that if
$$
w \in W, \ \gamma(w) < \infty,  \eqno(1.12)
$$
then

$$
1 \cdot||f||^*_w  \le ||f||_w \le \gamma(w) \cdot ||f||^*_w, \eqno(1.13)
$$
 and both  the coefficients $ "1" $ and $ "\gamma(w)" $ in (1.13) are the best possible. \par

\vspace{3mm}

 On the other word, the space $  Y = Y_w, \ w \in W $ in this example consists on all the right continuous  functions
defined on the set $ [0, \mu(X) ] $  equipped with the norm

$$
||g||Y_w = \sup_t [|g(t)| w(t)]
$$
and moreover the exact values of constants $ K_1(Y_w), \  K_2(V_w) $   are correspondingly:
$  K_1(Y_w) = 1, \ K_2(Y_w) = \gamma(w).  $\par

\vspace{4mm}

 We use the symbols $C(X,Y),$ $C(p,q;\psi),$ etc., to denote positive finite
constants along with parameters they depend on, or at least
dependence on which is essential in our study. To distinguish
between two different constants depending on the same parameters
we will additionally enumerate them, like $C_1(X,Y)$ and $C_2(X,Y).$

\vspace{4mm}
{\bf Layout of the paper.}
\vspace{3mm}
 The paper is organized as follows. In the next section we study the  estimations of Hardy operators
in weighted Lebesgue-Riesz spaces.   Third section is devoted to the  multidimensional  generalization
of  Hardy operators  estimations. \par
 The fourth section contains the main result of offered article:  sufficient conditions for normability
of generalized Lorentz spaces. In the fifth section we investigate the boundedness of Hardy operator in the
so-called anisotropic Grand Lebesgue spaces. \par
 In addition,  the last section contains a few review about properties of offered here spaces,
 in particular, calculation its fundamental function. \par

\hfill $\Box$ \\
 \bigskip

\section{Auxiliary facts: estimations of Hardy operators.}

\vspace{3mm}

  Let again $ (V, ||\cdot||V ) $ be the Banach functional space defined on the set $ R^+ = (0,\infty).  $
 Recall that the classical Hardy's operator  $ H = H[f] = H[f](t) $ (on the other term, Cesaro average) is
 defined as follows:

 $$
 H[f](t) = \frac{1}{t} \ \int_0^t f(s) ds. \eqno(2.1)
 $$
 {\it  It will be presumed that the Hardy's operator is defined on the space  } $ V $ {\it  and is bounded therein:}

$$
||H||_{V \to V} \stackrel{def}{=} \sup_{f \in V, ||f||V = 1  } ||H[f]||V < \infty. \eqno(2.2)
$$

 For instance, if $ V = L_p(R^+), \ 1 < p  \le \infty, $ then   $ ||H||_{L_p(R^+) \to L_p(R^+)} \le p':= p/(p-1),   $ and
the last estimation is not improvable  \cite{Hardy1}. \par
 There are many estimations for the norm of Hardy's operators in different spaces, for instance, in weight Lebesgue spaces
$ L_p(b): $

$$
||f||L_p(b) := \left[ \int_0^{\infty} |f(t)|^p \ b(t) \ dt  \right]^{1/p}, \eqno(2.3)
$$
where $ b = b(t)  $ is non-negative measurable local integrable function (weight), see for example
\cite{Arino1},  \cite{Bennet2}, \cite{Bloom1}, \cite{Bradley1}, \cite{Edmund1}, \cite{Hardy1},
\cite{Kufner2}, \cite{Maligranda11}, \cite{Maz'ja1}, \cite{Martin1}, \cite{Mitrinovic1},
\cite{Muckenhoupt1}, \cite{Stepanov1}, \cite{Xiao1} etc.; in Grand Lebesgue spaces \cite{Fiorenza3}, \cite{Ostrovsky112}; in
generalized weight Lorentz spaces \cite{Carro22}, \cite{Talenti1}, \cite{Tomaselli1} etc.  \par

 Note that in the article  of S.Bloom and R.Kerman  \cite{Bloom1} is considered more general  operator of a view, e.g.

 $$
 I_{\alpha, \beta}[f](x) = x^{-\beta} \int_0^x (x-y)^{\alpha} f(y) dy.
 $$

  P.R.Beesack in \cite{Beesack1} established the following result.  Let $ s=s(t) \ge 0, t > 0 $ be non-negative decreasing
function  such that

$$
S(x) := \int_0^x s(t) dt > 0, \ x > 0;
$$
$ f = f(t) $ be any function from the space $ L_p(1), \ p \in  (1,\infty). $ Then

$$
\left[ \left( \frac{1}{S(x)}  \int_0^{x} s(x-t) f(t) dt       \right)^p  \right]^{1/p} \le
\frac{p^2}{(p-1)} \left[ \int_0^{\infty} |f(t)|^p dt \right]^{1/p}.
$$

\vspace{3mm}

  {\it We refer here also the famous result of J.S. Bradley \cite{Bradley1}.} The inequality of a view
 with weights $ u^p(x), \ v^p(x) $
 $$
 \left( \int_0^{\infty} \left| \int_0^x f(s) \ ds  \right|^p  \ u^p(x) \ dx \right)^{1/p} \le C_p(u,v) \
 \left( \int_0^{\infty} |f(x)|^p  \ v^p(x) \ dx   \right)^{1/p} \eqno(2.4)
 $$
 is true for arbitrary (measurable) function $ f: R^+ \to R, $ where constant $ C_p(u,v)$ does not
depend on $ f, \ f \in L_p(v^p) $ iff

$$
B_p(u,v):= \sup_{r > 0} \left[ \left( \int_r^{\infty} u^p(x) dx \right)^{1/p} \cdot
\left( \int_0^r v^{-p'}(x) \ dx  \right)^{1/p'}  \right] < \infty \eqno(2.5)
$$
and moreover

$$
 B_p(u,v) \le C_p(u,v) \le p^{1/p} \ (p')^{1/p'} \ B_p(u,v).
$$
 Note that $ 1 \le p^{1/p} \ (p')^{1/p'} \le 2, $ therefore  $ B_p(u,v) \le C_p(u,v) \le 2B_p(u,v).  $\par
\vspace{3mm}
 For the different powers, i.e. for the inequality of a view

 $$
 \left( \int_0^{\infty} \left| \int_0^x f(s) \ ds  \right|^q  \ u^q(x) \ dx \right)^{1/q} \le C_{p,q}(u,v)) \
 \left( \int_0^{\infty} |f(x)|^p  \ v^p(x) \ dx   \right)^{1/p} \eqno(2.6)
 $$
 J.S. Bradley proved in \cite{Bradley1} that if  $ 1 < p \le q < \infty $ and

 $$
 B_{p \le q}(u,v) \stackrel{def}{=} \sup_{r>0} \left[ \left(\int_r^{\infty} u^q(x) \ dx \right)^{1/q} \cdot
  \left(\int_0^r v^{-p'}(s) \ ds  \right)^{ 1/p' }  \right]< \infty, \eqno(2.6a)
 $$
 then $ C_{p,q}(u,v) := C_{p \le q}(u,v) < \infty; $ moreover

 $$
 B_{p \le q}(u,v) \le C_{p \le q}(u,v) \le p^{1/q} \ (p')^{1/p'} B_{p \le q}(u,v).
 $$

 The case $  1 \le q < p < \infty  $  in the inequality (2.6) was investigated by  V.G.Maz'ja in
\cite{Maz'ja1}, chapter 11.  Indeed, the assertion (2.6) holds true under  condition $  1 \le q < p < \infty  $
for arbitrary admissible function  $ f: (0, \infty) \to R $ iff

$$
B_{p > q}(u,v):= \left\{\int_0^{\infty} \left[ \left( \int_0^x |v(y)|^{-p'} dy \right)^{q-1} \int_x^{\infty} |u(y)|^q dy \right]^{p/( p-q )  }
   \frac{dx}{|v(x)|^{p'}}  \right\}^{(p-q)/(pq)} < \infty \eqno(2.7)
$$
herewith

 $$
 \left[\frac{p-q}{p-1} \right]^{(q-1)/q} \ q^{1/q} \ B_{p > q}(u,v) \le C_{p >q}(u,v)) \le
 \left[ \frac{p}{p-1} \right]^{(q-1)/q} \ q^{1/q} \ B_{p > q}(u,v).
\eqno(2.7a)
$$

\vspace{3mm}
 We refer also the following important for us result belonging to V.D.Stepanov  \cite{Stepanov1} relative the
{\it non-increasing non-negative } function $ f. $ Consider the inequality of a view

$$
||H[f]||L_q(w) \le C(w,v) \cdot ||f||L_p(v), \ 0 < C(w,v) < \infty,
$$
and define $ p' = p/(p-1), \ p > 1; \ V(x) = \int_0^x v(s) ds, \ 1/r = 1/q - 1/p,$

$$
A_0 = \sup_{t > 0} \left[ \left( \int_0^t w(x) \ dx  \right)^{1/q} \cdot \left( \int_0^t v(x) \ dx \right)^{-1/p}      \right],
$$

$$
A_1 = \sup_{t > 0} \left[ \left(\int_t^{\infty} x^{-q} \ w(x) \ dx  \right)^{1/q} \cdot
 \left( \int_0^t x^{p'} \ V^{-p'}(x) \ v(x) \ dx   \right)^{1/p'}   \right],
$$

$$
B_0=\left\{ \int_0^{\infty} \left[  \left( \int_0^t w(x) dx   \right)^{1/p} \cdot
\left(\int_0^t v(x) dx \right)^{-1/p} \right]^r w(t) dt  \right\}^{1/r}, \eqno(2.8)
$$

$$
B_1=\left\{ \int_0^{\infty} \left[  \left( \int_t^{\infty}x^{-q} w(x) dx   \right)^{1/q} \cdot
\left(\int_0^t x^{p'} \ V^{-p'}(x) v(x) dx \right)^{1/q'} \right]^r  \ t^{p'} V^{-p'}(t) v(t) \ dt  \right\}^{1/p}.\eqno(2.9)
$$

 If  $ 1 < p \le q < \infty, $ then

 $$
 \alpha_1(p,q) (A_0 + A_1) \le C(w,v) \le \alpha_2(p,q) (A_0 + A_1), \ 0 < \alpha_1(p,q) \le \alpha_2(p,q) < \infty;\eqno(2.10)
 $$
if $ 1 < q < p < \infty, $ then

$$
\beta_1(p,q)(B_0 + B_1) \le    C(w,v) \le \beta_2(p,q)(B_0 + B_1), \ 0 < \beta_1(p,q) \le \beta_2(p,q) < \infty.\eqno(2.11)
$$
  The "equal" case when $ p=q $ and $ u = v $ was considered, e.g. in  \cite{Arino1}, \cite{Bennet2}, \cite{Muckenhoupt1}.
Namely, the inequality

$$
||H[f]||L_p(w) \le C_p(w) \cdot ||f||L_p(w) \eqno(2.12)
$$
holds true iff

$$
 \exists D = D(p,w) \in (0,\infty), \ \int_0^{\infty} s^{-p} w(s) ds \le D(p,w) \ t^{-p} \int_0^t w(s) \ ds. \eqno(2.13)
$$

 See also \cite{Arino1}, \cite{Bennet2}, \cite{Muckenhoupt1} where are obtained alike result  without constants computation.
 In the theses  of L.Arendarenko \cite{Arendarenko1} and O.Popova \cite{Popova1}  there is a comprehensive review
 about this problem and are offered some new results. \par

\vspace{3mm}

{\bf Example 2.1.}\par
 Let us consider an inequality of a view:

 $$
 ||x^{\alpha} \ H[f]||_q \le K_{\alpha, \beta}(p) \ ||x^{\beta} \ f||_p, \ \alpha,\beta = \const, \eqno(2.14)
 $$
 or equally

 $$
 || \ H_{\alpha, \beta}[g] \ ||_q \le K_{\alpha, \beta}(p) \ || \ g \ ||_p, \ \alpha,\beta = \const, \eqno(2.14a)
 $$
 where

 $$
 H_{\alpha, \beta}[g](x) \stackrel{def}{=} \ x^{\alpha} \ H[x^{- \beta} \ g](x).
 $$

In the capacity of  the value $ K_{\alpha, \beta}(p) $ we understood its minimal (and  implied to be finite) value:

$$
K_{\alpha, \beta}(p) = \sup_{ 0 < ||x^{\beta} f||_p < \infty } \frac{||x^{\alpha} \ H[f] \ ||_q }{||x^{\beta}  f||_p }  =
 \sup_{ 0 < ||g||_p < \infty } \frac{|| \ H_{\alpha,\beta}[g] \ ||_q }{|| \  g||_p}.\eqno(2.15)
$$

\vspace{3mm}

{\bf A.} We investigate first of all the case $ \alpha \ge \beta $  or equally $ q \ge p $  (case "Bradley").\par
  Denote  for simplicity

$$
p_0 = 1/(1-\beta), \ p_+ = 1/(\alpha-\beta),  \ q_0 = 1/(1-\alpha), \ q_+ = + \infty, \ \delta = \alpha - \beta.
$$
 It is clear that $ 0 < \delta \le 1 $ and $  p_0 < p \le p_+ \ \Leftrightarrow q_0 < q \le \infty.  $\par

  It follows from Bradley's inequality then (2.14) holds iff

$$
0  \le \alpha, \beta < 1, \  \delta = \alpha - \beta =  \frac{1}{p} - \frac{1}{q}, \eqno(2.16a)
$$

$$
p > \frac{1}{1-\beta} = p_0, \ q > \frac{1}{1-\alpha} = q_0. \eqno(2.16b)
$$

 Since $ q < \infty, $ we conclude in the considered case $ p \le 1/(\alpha - \beta) = p_+. $ So,

$$
p_0 < p \le p_+, \hspace{5mm} q_0 < q \le q_+.
$$

\vspace{3mm}
 {\bf Remark 2.1.}  Note that the values $ \alpha,\beta, p_0, p_+, q_0, q_+ $ presumed to be constants, but the values
 $ p,q $ are variable.\par
\vspace{3mm}

 We deduce from the Bradley's inequality:

$$
K_{\alpha, \beta}(p) \le C(\alpha,\beta)  \left[ \frac{p}{p-p_0} \right] ^{1 - \alpha + \beta}, \ p_0 < p \le p_+.\eqno(2.16c)
$$

 As a particular case: as  $  \alpha = \beta $

$$
K_{\alpha, \alpha}(p) \le  \frac{C_2(\alpha) \ p }{p-p_0}, \ p_0 < p \le p_+.\eqno(2.16d)
$$
 We will prove further the inverse  inequality. Thus:

$$
K_{\alpha, \alpha}(p) \asymp  \frac{p}{ p-p_0}, \ p_0 < p \le p_+.
$$

 Note that in the case $ \alpha = \beta = 0 $ we obtain the weak version of the classical Hardy's inequality with coefficient

  $$
K_{0, 0}(p) \asymp \frac{p}{p-1}, \ p \in (1,\infty],
$$
 which is exact up to multiplicative constant.\par

\vspace{3mm}
{\bf B.} Case "Maz'ja".  Let now  $ \alpha < \beta $  or equally $ q < p. $  \par
\vspace{3mm}
 Note that in the considered here restriction  and when $ u(x) := u_{\alpha}(x) = x^{\alpha}, \ v(x) := v_{\beta}(x) = x^{\beta}  $

 $$
 B_{p > q} \left(u_{\alpha}(\cdot),v_{\beta}(\cdot) \right) = +\infty.
 $$
 Therefore, the inequality (2.6) with $ u(x) = u_{\alpha}(x) = x^{\alpha}, \ v(x) = v_{\beta}(x) = x^{\beta} $
  or  equally $ C_{q > p}\left(u_{\alpha}(\cdot),v_{\beta}(\cdot) \right)  < \infty $ may be true iff  $ \alpha > \beta $  or equally $ q > p. $  \par
 Outcome:
\vspace{3mm}

{\bf Proposition 2.1.} {\it The constant} $ K_{\alpha, \beta}(p) = K_{\alpha, \beta}(p,q) $ {\it from inequality (2.15) is finite  iff}
$$
 0 \le \alpha, \beta < 1, \hspace{5mm}  \alpha > \beta, \eqno(2.17a)
$$

$$
 \alpha - \beta =  \frac{1}{p} - \frac{1}{q}, \eqno(2.17b)
$$

$$
p_0 = \frac{1}{1-\beta} < p \le \frac{1}{\alpha-\beta} = p_+, \  q_0 = \frac{1}{1-\alpha} < q \le \infty = q_+.\eqno(2.17c)
$$

\vspace{3mm}

 More general weight inequality with exact value of the constant  belongs to G.Hardy (\cite{Bennet1}, p. 124-125):

 $$
  \left\{\int_0^{\infty} \left(t^{\nu} H[f](t) \right)^q \ dt/t \right\}^{1/q} \le \frac{1}{1-\nu} \
  \left\{\int_0^{\infty} (t^{\nu} \ f(t) )^q  dt/t  \right\}^{1/q}, \ \nu = \const < 1.
 $$

\vspace{3mm}

  The accuracy calculation used Bradley's  estimation tell us that $  K_{\alpha, \beta}(p) \ge K^0_{\alpha, \beta}(p), $
where

$$
K^0_{\alpha, \beta}(p)= (1-\beta)^{\delta -1} \cdot (p-1)^{1-1/p} \cdot \delta^{1/p - \delta} \cdot
 \frac{[p_+ - p]^{1/p - \delta}}{[p-p_0]^{1-\delta}} =
$$

 $$
 (1-\beta)^{\alpha - \beta - 1} \cdot (p-1)^{1-1/p} \cdot
  [p - 1/(1-\beta)]^{\alpha - \beta - 1}  \cdot [1 + p(\beta-\alpha)]^{\alpha - \beta - 1/p},
 $$
and

$$
K_{\alpha, \beta}(p) \le  p^{1/q} \ (p')^{1/p'} \cdot K^0_{\alpha, \beta}(p), \ q = q(p) = p/(1 - p\delta ).
$$
  As long as

 $$
 p^{1/q} \ (p')^{1/p'} = (p')^{1/p'} \cdot p^{1/p + \beta - \alpha}  \le \left[(p')^{1/p'} \cdot p^{1/p} \right] \cdot
  p^{ \beta - \alpha} \le  2,
 $$
 we conclude

  $$
 K_{\alpha, \beta}(p) \le  K^+_{\alpha, \beta}(p) \stackrel{def}{=} 2 \  K^0_{\alpha, \beta}(p).
  $$
  \vspace{3mm}

  {\it  Our hypothesis:  } under our restrictions (2.17a),  (2.17b), (2.17c)
    $$
 K_{\alpha, \beta}(p) = C(\alpha,\beta) \cdot K^0_{\alpha, \beta}(p), \   1 \le C(\alpha,\beta) \le 2.
  $$
    \vspace{3mm}
  See also  \cite{Okikiolu1}, pp. 211-221. \par
\vspace{3mm}

{\bf Example 2.2.}\par
 More generally, consider the inequality of a view

 $$
 ||x^{\alpha} \ L(x) \ H[f]||_q \le K_{L,M; \alpha, \beta}(p) \ ||x^{\beta} \ M(x) \ f||_p, \ \alpha,\beta = \const,
\eqno(2.18)
 $$
where $  L(x), M(x) $ are slowly varying simultaneously as $ x \to 0+ $ and as $ x \to \infty $ continuous in the
semi-axis $ (0, \infty) $ positive functions.  As ordinary, in the capacity of  the value $ K_{L,M; \alpha, \beta}(p) $ we understood
its minimal value, presumed to be finite.\par
 We conclude using at the same method as in last example  that the estimate (2.18) holds true iff $ (\alpha, \beta, p,q) $
satisfy  the conditions   (2.17a), (2.17b), (2.17c) and

$$
0 < \inf_{x > 0} \left[\frac{L(x)}{M(x)} \right] \le \sup_{x > 0} \left[\frac{L(x)}{M(x)} \right] < \infty. \eqno(2.18a)
$$

 Moreover, under these conditions
$$
K_{L,M;\alpha, \beta}(p) \asymp C(L,M; \alpha,\beta) \cdot \left\{ \frac{p}{ [p-1/(1-\beta)]} \right\}^{1 - \alpha + \beta}. \eqno(2.18b)
$$

 We use the {\it multidimensional} generalization of the so-called {\it dilation }, or {\it scaling} method, see  \cite{Talenti0},
 \cite{Stein1}, chapter 3.  Indeed, let us introduce the following dilation operator $ T_{\lambda}[f](x) := f(\lambda x), \lambda \in (0,\infty).  $
 Suppose the inequality (2.18) is true for arbitrary function $ f $ from the Schwartz space $ S(0,\infty), $   and substitute in
(2.18) the function $ T_{\lambda}[f] \in S(0,\infty) $ instead $ f, \ f \ne 0: $

$$
 ||x^{\alpha} \ L(x) \ H[T_{\lambda}[f]]||_q \le K_{L,M; \alpha, \beta}(p) \ ||x^{\beta} \ M(x) \ T_{\lambda}[f]||_p.\eqno(2.18c)
$$

We deduce consequently:

$$
||x^{\beta} \ M(x) \ T_{\lambda}[f]||_p^p= \int_0^{\infty} x^{\beta p} M^p(x) |f(\lambda x)|^p dx =
$$

$$
\lambda^{-1-\beta p} \int_0^{\infty} y^{\beta p} M^p(y/\lambda) |f(y)|^p dy \asymp \lambda^{-1-\beta p} M^p(1/\lambda)
\int_0^{\infty} y^{\beta p} |f(y)|^p dy;
$$

$$
||x^{\beta} \ M(x) \ T_{\lambda}[f]||_p \asymp \lambda^{-1/p - \beta} M(1/\lambda) || \ x^{\beta} f \ ||_p;
$$

$$
H [T_{\lambda}[ f]] =  T_{\lambda}[ H [f]];
$$

$$
 ||x^{\alpha} \ L(x) \ H[T_{\lambda}[f]]||_q \asymp \lambda^{-1/q - \alpha} L(1/\lambda) || \ x^{\alpha} H[f] \ ||_q;
$$

$$
 \lambda^{-1/q - \alpha} L(1/\lambda) || \ x^{\alpha} H[f] \ ||_q \asymp \lambda^{-1/p - \beta} M(1/\lambda) || \ x^{\beta} f \ ||_p;
$$
therefore

$$
1/q + \alpha = 1/p + \beta, \  L(1/\lambda)  \asymp M(1/\lambda),
$$
which is equivalent to our assertion.\par
 The passing to the limit as $ \lambda \to 0+ $ or $ \lambda \to \infty $ in the considered case is grounded in
 \cite{Liflyand1}, \cite{Ostrovsky2}. \par
 The relation (2.18b) follows immediately from Bradley's estimation by means if properties of slowly and regular
varying functions,  see \cite{Bingham1}, chapter 3; \cite{Seneta1}, chapter 2. \par

\vspace{3mm}

{\bf Example 2.3.}\par
 The {\it lower} bound in the inequality (2.16c) may be obtained even without restriction $ \beta < \alpha $
 by means of consideration of an example

$$
f_0(x) = x^{-1} \ (\log x)^{\Delta} \ I_{(1,\infty)}(x).
\eqno(2.19)
$$
  Here and further $ I_A(x) = 1, \ x \in A,  \ I_A(x) = 0, \ x \notin A.  $ \par
 Namely, it is easy to compute that under formulated before conditions and restrictions and as $ \Delta = \const >> 1 $
 there holds

$$
 \frac{||x^{\alpha} \ H[f_0]||_q }{||x^{\beta} \ f_0||_p } \asymp \left[\frac{p}{p-1/(1-\beta)} \right]^{1 - \alpha + \beta},
 \ p > 1/(1-\beta).
\eqno(2.19)
$$

 In detail,

 $$
 || \ x^{\beta} \ f_0 \ ||_p^p = \int_1^{\infty} x^{\beta p - p} \ \log^{\Delta p} x \ dx =
 \frac{\Gamma(\Delta p + 1)}{[p(1-\beta) - 1 ]^{\Delta + 1/p}}; \eqno(2.20)
 $$

$$
x^{\alpha} H[f_0](x) = x^{\alpha -1} \ \int_1^x  s^{-1} \ \log^{\Delta} s \ ds \ I_{(1,\infty)}(x) =
$$

$$
(1+\Delta)^{-1} \ I_{(1,\infty)}(x) \ x^{\alpha-1} \ \log^{\Delta+1}x; \eqno(2.21)
$$

$$
||x^{\alpha} H[f_0]||_q = \frac{\Gamma^{1/q}((\Delta+1)q+1)}{(\Delta+1)[q(1-\alpha) - 1 ]^{\Delta + 1 + 1/q} }; \eqno(2.22)
$$

$$
\lim_{\Delta \to \infty} \frac{||x^{\alpha} \ H[f_0]||_q }{||x^{\beta} \ f_0||_p }= \left[\frac{1}{p-1/(1-\beta)} \right]^{1 - \alpha + \beta} \
\cdot (1-\beta)^{\alpha - \beta - 1} \cdot \frac{p}{[1-p(\alpha - \beta)]^{\alpha - \beta - 1/p}}.\eqno(2.23)
$$

 We used the Stirling's formula for the Gamma function $ \Gamma(z), \ z \to \infty. $ \par
 \vspace{3mm}

 {\bf Corollary 2.1.} \\
  $$
 K_{\alpha,\beta}(p) \ge \left[\frac{1}{p-1/(1-\beta)} \right]^{1 - \alpha + \beta} \
  \cdot (1-\beta)^{\alpha - \beta - 1} \cdot \frac{p}{[1-p( \alpha - \beta)]^{\alpha - \beta - 1/p}}.\eqno(2.24)
 $$

  \vspace{3mm}
  {\bf Corollary 2.2.} \\
 As long as $ K^{0}_{\alpha,\beta}(p) $  is less than the right-hand side of inequality (2.24), we conclude that
in general case the function   $ K^{0}_{\alpha,\beta}(p) $  is not exact lower bound in the Bradley bilateral inequality. \\

\hfill $\Box$ \\
 \bigskip

\section{Multidimensional case.}

\vspace{3mm}

 We consider further in this section the so-called $ d - $ dimensional Hardy's operator $ H_d[f] $ defined on the functions
defined on the "octant"   $  R^d_+ = (R^1_+)^d $ by a formula

$$
H_d[f](x_1,x_2, \ldots, x_d) = \frac{1}{x_1 \ x_2, \ \ldots, \ x_d} \cdot
\int_0^{x_1}\int_0^{x_2} \ldots    \int_0^{x_d} f(y_1,y_2, \ldots, y_d) \ dy_1  dy_2 \ldots  dy_d. \eqno(3.1)
$$
 The $ L_p(b) $ estimations for the norm   of $ H_d[\cdot] $ see in  \cite{Sawyer1}, \cite{Sedov1},
 \cite{Sinnamon1}. \par
  Another approach ("spherical definition") see in \cite{Christ1}, \cite{Faris1}. \par
 In \cite{Faris1} are described in addition an applications of these estimations into the quantum mechanic.\par
\vspace{3mm}

 We intent to investigate in  this section the inequality of a view

 $$
 || w(x) \ H_d[f](x) ||_q \le C(w, v) \cdot || v(x) \ f(x) ||_p, \   0 < C(w,v) < \infty,  \eqno(3.2)
 $$
where the weight function $ w(x) = w_{\alpha}(x) $ is homogeneous of degree $ \alpha $ continuous on the unit sphere positive function,
 the weight function $ v(x) = v_{\beta}(x) $ is homogeneous of degree $ \beta $ continuous on the unit sphere positive function.\par

\vspace{3mm}

{\bf Proposition 3.1.}{\it Suppose the inequality (2.15) holds true  for each non-zero
 function   from the Schwartz class $ S(R^d_+): \ f \in S(R^d_+). $ Then}

$$
\alpha - \beta = d \left(\frac{1}{p} - \frac{1}{q} \right). \eqno(3.3)
$$
{\bf Proof.}   We will use again the so-called  {\it dilation,} or equally {\it scaling} method, see \cite{Stein1}, chapter 10;
\cite{Talenti0}.  Namely, let us define the {\it family} of dilation operators $ T_{\lambda} = T_{\lambda}[f] $ as follows:

$$
T_{\lambda}[f](x) = f(\lambda x), \ \lambda \in (0,\infty). \eqno(3.4)
$$
 Evidently, $  f_{\lambda}: = T_{\lambda}[f] \in S(R^d_+). $  We have:

 $$
 || \ w_{\alpha}(x) \ H_d[f_{\lambda}](x) \ ||_q \le || \ v_{\beta}(x) \ f(\lambda x),  \  ||_p. \eqno(3.5)
 $$

   Note that

 $$
 || \ v_{\beta}(x) \ f(\lambda x)  \  ||_p^p = \int_0^{\infty} v^p_{\beta}(x) \ |f(\lambda x)|^p \ dx =
 $$

$$
\lambda^{ -\beta p - d } \int_0^{\infty} v^p_{\beta}(y) \ |f(y)|^p \ dy = \lambda^{ -\beta p - d } \ || v_{\beta}(x) \ f(x) \ ||_p^p, \eqno(3.6)
$$
therefore

$$
|| \ v_{\beta}(x) \ f(\lambda x)  \  ||_p =
 \lambda^{ -\beta  - d/p } \ || v_{\beta}(x) \ f(x) \ ||_p. \eqno(3.7)
$$
 Further, $ H_d[T_{\lambda }f] = T_{\lambda }H_d[f],   $

$$
 || w_{\alpha}(x) \ H_d[T_{\lambda}[f]] ||^q_q = || w_{\alpha}(x) \ T_{\lambda}[H_d[f]] ||^q_q=
\lambda^{-\alpha q -d} || \ w_{\alpha}(x) \ H_d[f_{\lambda}](x) ||_q^q,
$$

$$
 || w_{\alpha}(x) \ H_d[T_{\lambda}[f]] ||_q = \lambda^{-\alpha  -d/q} || \ w_{\alpha}(x) \ H_d[f](x) ||_q. \eqno(3.8)
$$

We get substituting into inequality (3.4):

$$
\lambda^{-\alpha  -d/q} || \ w_{\alpha}(x) \ H_d[f](x) ||_q \le  C(w_{\alpha},v_{\beta}) \cdot
\lambda^{ -\beta  - d/p } \ || v_{\beta}(x) \ f(x) \ ||_p.\eqno(3.9)
$$
 Since the value $  \lambda $ is arbitrary positive, we conclude from (3.9)

$$
\alpha + d/q = \beta + d/p,
$$
which is equivalent to (3.3). \par
 Now we investigate the inequality (2.14) without assumption  of homogeneity of a functions $ w(x), v(x). $  Namely,
we suppose the existence of finite constants $ \alpha(0), \alpha(\infty), \beta(0), \beta(\infty)  $  for which the following functions

$$
v_0(x) := \sup_{\lambda \in (0,1)} \frac{v(\lambda x)}{\lambda^{\beta(0)}}, \
v_{\infty}(x) := \inf_{\lambda \in (1, \infty)} \frac{v(\lambda x)}{\lambda^{\beta(\infty)}},
$$

$$
w_0(x) := \inf_{\lambda \in (0,1)} \frac{w(\lambda x)}{\lambda^{\alpha(0)}}, \
w_{\infty}(x) := \sup_{\lambda \in (1, \infty)} \frac{w(\lambda x)}{\lambda^{\alpha(\infty)}},
$$
are non-trivial: non-zero and integrable.\par

\vspace{3mm}

{\bf Proposition 3.2.}{\it Suppose in addition the inequality (3.2) holds true  for each non-zero
 function   from the Schwartz class $ S(R^d_+): \ f \in S(R^d_+). $ Then}

$$
\alpha(0) - \beta(0) \ge d \left(\frac{1}{p} - \frac{1}{q} \right).
$$

$$
\alpha(\infty) - \beta(\infty) \le d \left(\frac{1}{p} - \frac{1}{q} \right).
$$
{\bf Proof} is at the same as in proposition (2.1) and may be omitted. \par
 Obviously, when $ \alpha(0) - \beta(0) = \alpha(\infty) - \beta(\infty),  $ then

$$
\alpha(0) - \beta(0) = \alpha(\infty) - \beta(\infty) = d \left(\frac{1}{p} - \frac{1}{q} \right).
$$

\vspace{3mm}
 We recall here the definition of
the so-called anisotropic Lebesgue (Lebesgue-Riesz) spaces, or equally the spaces with mixed norms. More
detail information about this spaces see in the books  of Besov O.V., Il’in V.P., Nikol’skii S.M.
\cite{Besov1}, chapter 16,17; Leoni G. \cite{Leoni1}, chapter 11; using for us theory of
operators interpolation in this spaces see in \cite{Besov1}, chapter 17,18. \par

  Let $ (X_j,A_j,\mu_j), \ j=1,2,\ldots,d $ be measurable spaces with sigma-finite
non - trivial measures $ \mu_j. $
Let also $ p = (p_1, p_2, . . . , p_d) $ be $ d- $ dimensional vector such that
$ 1 \le p_j \le \infty.$ \par

 Recall that the anisotropic Lebesgue space $ L_{ \vec{p}} $ consists on all the  total measurable
real valued function  $  \ f = f(x_1,x_2,\ldots, x_d) = f( \vec{x} ), $

$$
f:  \otimes_{j=1}^d X_j \to R
$$

with finite norm $ \ |f|_{ \vec{p} } \stackrel{def}{=} $

$$
\left( \int_{X_d} \mu_d(dx_d) \left( \int_{X_{d-1}} \mu_{d-1}(dx_{d-1}) \ldots \left( \int_{X_1}
 |f(\vec{x})|^{p_1} \mu_1(dx_1) \right)^{p_2/p_1 }  \ \right)^{p_3/p_2} \ldots   \right)^{1/p_d}. \eqno(3.10)
$$

 Note that in general case $ |f|_{p_1,p_2} \ne |f|_{p_2,p_1}, $
but $ |f|_{p,p} = |f|_p. $ \par

 Observe also that if $ f(x_1, x_2) = g_1(x_1) \cdot g_2(x_2) $ (condition of factorization), then
$ |f|_{p_1,p_2} = |g_1|_{p_1} \cdot |g_2|_{p_2}, $ (formula of factorization). \par

  We use here the case $ X_j = R_+, \ \mu_j(dx_j) = dx_j. $\par

 \vspace{3mm}

 We consider in this section the weight multidimensional (vector): $ d \ge 2 $ generalization of
weight Hardy's $ L_p(b_1) \to L_q(b_2) $ estimations.\par
 In this section $ x = \vec{x} \in R^d $ be $ d- $ dimensional vector, $ d =  2,3, . . .  $ which consists
on the $ d $  coordinates  $ {x_j}, j = 1, 2, . . . , d : \ $
$$
x = (x_1, x_2, . . . , x_d),
$$
$$
 \alpha = \vec{\alpha} = \{\alpha_1, \alpha_2, \ldots, \alpha_d \},
 \beta = \vec{\beta} = \{\beta_1,\beta_2, \ldots,\beta_d \}.
$$

 We denote as ordinary
$$
x^{\alpha} = \vec{x}^{\vec{\alpha}} = \prod_{j=1}^d x_j^{\alpha_j}, \hspace{5mm}
y^{\beta} =  \vec{y}^{\vec{\beta}} = \prod_{j=1}^d y_j^{\beta_j}.
$$

\vspace{3mm}

  Let $ f, \ f: R^d \to R $ be (total) measurable function.
 Let also
 $$
  \alpha_j,\beta_j = \const \in [0,1), \ \alpha_j > \beta_j,
 \ j=1,2,\ldots,d; \ p_j \in ( (1/(1-\beta_j), 1/(\alpha_j - \beta_j)). \eqno(3.11)
 $$

We define the function $ q_j=q_j(p_j) $ as follows:
$$
 \frac{1}{p_j} - \frac{1}{q_j}  = \alpha_j-\beta_j.  \eqno(3.11a)
$$

 The equality (3.11a) defines the dependance between $ \vec{p}  $ and $  \vec{q}; $  we will denote
this functions as follows:

$$
\vec{p}=  \vec{p}( \vec{q} ), \ \vec{q}=  \vec{q}( \vec{p} ).
$$
 Obviously,  two functions $ \vec{p}=  \vec{p}( \vec{q} ) $ and $  \vec{q}=  \vec{q}( \vec{p} ) $ are
reciprocal inverse.   \par

\vspace{3mm}

{\bf Theorem 3.2.} \par
 {\it The conditions (2.17a), (2.17b), (2.17c) for  the variables $ (\alpha_j, \beta_j, p_l, q_j) $
  are necessary and sufficient for the existence of non-trivial coefficient}
$ K(d; \vec{\alpha}, \vec{\beta},\vec{p}) $ {\it for the following estimate:}

 $$
 || \vec{x}^{\vec{\alpha}} \ H_d[f](x) ||_{\vec{q}} \le K(d;\vec{\alpha},\vec{\beta};\vec{p} ) \cdot
 || \vec{x}^{\vec{\beta}} \ f(x) ||_{\vec{p}}, \   0 < K(d,\vec{\alpha},\vec{\beta} ) < \infty,  \eqno(3.12)
 $$
{\it and under this conditions  for the minimal value of coefficient}
$  K(d; \vec{\alpha}, \vec{\beta}; \vec{p}) $ {\it there hold the
following equality:}

$$
 K(d; \vec{\alpha}, \vec{\beta}; \vec{p})  =  \prod_{j=1}^d K_{\alpha_j,\beta_j}(p_j), \eqno(3.13a)
$$

 {\it and } $ K(d; \vec{\alpha}, \vec{\beta},\vec{p}) = \infty  $ {\it in other case.} \par

\vspace{3mm}

{\bf The proof} is at the same as one in \cite{Ostrovsky114} for weight Riesz and Fourier integral transform.
 Namely, let us introduce the following one-dimensional operators

 $$
 H^{(j)}[f](x_1,x_2, \ldots, x_{j-1}, x_j, x_{j+1}, \ldots, x_d) = \frac{1}{x_j} \int_0^{x_j}
 f(x_1,x_2, \ldots, x_{j-1}, s_j, x_{j+1}, \ldots, x_d) \ ds_j, \eqno(3.14)
 $$
then

$$
H_d = \otimes_{j=1}^d H^{(j)}. \eqno(3.15)
$$
 It is sufficient to consider only the two-dimensional case $ d = 2. $ Denote

 $$
 z(x_1,x_2) = H_2[f](x_1,x_2) = \frac{1}{x_1} \int_0^{x_1} g(s_1; x_2) d s_1 = H^{(1)}[g](x_1; x_2), \eqno(3.16)
 $$
where

$$
g(s_1; x_2) = H^{(2)}[f](s_1,x_2) = \frac{1}{x_2}\int_0^{x_2} f(s_1, x_2) \ d s_1. \eqno(3.17)
$$

  We have using the one-dimensional estimate (2.16c) for $ u(\cdot, \cdot): $

$$
|| x_1^{\alpha_1} \ z(\cdot, x_2)||_{q_1, X_1} \le  K_{\alpha_1, \beta_1}(p_1) \cdot
|| x_1^{\beta_1} \ g(\cdot, x_2)||_{p_1, X_1. } \eqno(3.18)
$$

 Now we apply the triangle inequality for the $ L(q_2) $ norm  and the one-dimensional estimate (2.16c):

$$
|| x_1^{\alpha_1} \ x_2^{\alpha_2} \ z(\cdot,\cdot)||_{q_1,q_2} = || x_1^{\alpha_1} \ || x_2^{\alpha_2} u(\cdot, x_2)||_{q_1, X_1} ||_{q_2, X_2}  \le
$$

$$
K_{\alpha_1, \beta_1}(p_1) \cdot  K_{\alpha_2, \beta_2}(p_2) \cdot
|| \  x_1^{\beta_1} \ || \ x_2^{\beta_2} f(\cdot, x_2)||_{p_1, X_1} \  ||_{p_2,X_2} =
$$

$$
K_{\alpha_1, \beta_1}(p_1) \cdot  K_{\alpha_2, \beta_2}(p_2) \cdot || \  x_1^{\beta_1} \ x_2^{\beta_2} \  f||_{p_1,p_2}. \eqno(3.19)
$$

 The {\it lower} estimate in (3.13)  may be obtained after consideration an example $ f_0(\vec{x}) $ of factorized function of a view

 $$
 f_0(\vec{x})  = \prod_{j=1}^d h_j(x_j). \eqno(3.20)
 $$

\vspace{3mm}

 It nay be considered analogously a more general case  of inequality of a view

 $$
 || \ u(\vec{x}) \ H_d[f] \ ||_{\vec{q}} \le C^{(d)}(u, v; \vec{p}, \vec{q}) \ || v(\vec{x}) \ f \ ||_{\vec{p}},  \eqno(3.21)
 $$
where both the weight $ u(\vec{x}) $ and $ v(\vec{x}) $ are {\it factorable:}

$$
u(\vec{x}) = \prod_{j=1}^d u_j(x_j), \hspace{5mm} v(\vec{x}) = \prod_{j=1}^d v_j(x_j).  \eqno(3.22)
$$

\vspace{3mm}

{\bf Theorem 3.3.} {\it  Let the weights functions } $ u(\vec{x}), \ v(\vec{x}) $ {\it  be factorable in the sense (3.22).
The inequality (3.21) holds true with the coefficient }

$$
C^{(d)}(u, v; \vec{p}, \vec{q}) = \prod_{j=1}^d C_{p_j,q_j}(u_j, v_j), \eqno(3.23)
$$
{\it where the function} $ C_{p,q}(u,v) $ {\it is defined in (2.6).}

\vspace{3mm}

{\bf Example 3.1.} Suppose

$$
u_j(x_j) = x_j^{\alpha_j} \ L_j(x_j), \hspace{5mm} v_j(x_j) = x_j^{\beta_j} \ M_j(x_j), \ j=1,2,\ldots, d, \eqno(3.24)
$$
where $ L_j(y), \  M_j(y), \ y \in (0, \infty)  $ are positive continuous slowly varying simultaneously as $ y \to 0+ $
and as $ y \to \infty $  functions.  The inequality (3.21) for these weights is valid iff the parameters
$ (\alpha_j, \beta_j; p_j, q_j) $ satisfy the conditions (3.11), (3.11a) and moreover

$$
0 < \min_j \inf_{y > 0} \left[\frac{u_j(y)}{v_j(y)} \right] \le \max_j \sup _{y > 0} \left[\frac{u_j(y)}{v_j(y)} \right] < \infty. \eqno(3.25)
$$

\hfill $\Box$ \\
 \bigskip

\section{Main result.}

\vspace{3mm}
 The following assertion is obvious. \par
\vspace{3mm}
{\bf Proposition 4.1.} The space $ (Y, ||\cdot||Y ) $ is the Banach (complete) rearrangement invariant functional space defined
on the set $ (X = \{x\}, \cal{A}, \mu).  $\par
\vspace{3mm}
 For instance, the triangle inequality  and homogeneity of the norm $ ||\cdot||V   $ follows immediately
from the equality (1.9). \par

  Denote as $ D_+ = D_+ (V)  $ the set of all positive decreasing (measurable) functions $ f: V \to R_+. $
 We define

 $$
 K(V) = K(V,H) \stackrel{def}{=} \sup_{0 \ne g \in D_+ (V)} \left[  \frac{||Hg||V}{||g||V} \right]. \eqno(4.1)
 $$
Obviously, $ K(V,H) \le ||H||_{(V \to V)}.  $

\vspace{3mm}
{\bf Theorem 4.1.}
$$
 1 \cdot ||| f |||^*_Y   \le ||| f |||_Y \le K(V,H) \cdot ||| f |||^*_Y, \eqno(4.2)
$$
{\it where both the constants "1" and $ " K(V,H)" $ are the best possible.  } \par

\vspace{3mm}
{\bf Proof.  Inequalities.} The left-hand side of proposition (3.2) follows immediately from the inequality
$ f^*(t) \le f^{**}(t). $ The main idea for proving of the right-hand side (3.2) is following:

$$
 f^{**}(t) = t^{-1} \int_0^t f^*(s) \ ds = H[f^*](t). \eqno(4.3)
$$
  Suppose $ K(V,H) < \infty; $ other case it is nothing to prove. \par
 Since the function $ t\to f^*(t) $ is positive and decreasing, we can use the definition of the constant
$ K(V,H): $

$$
||| f |||^*_Y = || f^{**}(\cdot) ||V \le K(V,H) \ || f^{*}(\cdot) ||V = K(V,H) \cdot ||| f |||^*_Y. \eqno(4.4)
$$

\vspace{3mm}
{\bf Proof.  Exactness.}   The  exactness of the left constant  "1"  is true, e.g.  for the spaces $ L_p, \ p \in (1,\infty), $
see \cite{Ostrovsky111}. But we can assert the exactness for each Banach functional space $ (V, ||\cdot||V), $ as in the article
\cite{Ostrovsky113}.  Indeed, let us denote

$$
\underline{K}(V) = \inf_{f \ne 0} \left[\frac{|||f|||V}{|||f|||^*V}   \right].
$$
 We introduce also the family of a functions of a view

 $$
f^*_h(t) = h_{\kappa}(t) = 1 - t^{\kappa}, \ t \in (0,1);
 $$
then

$$
f^{**}_h(t) = 1 - t^{\kappa}/(\kappa+1).
$$
 Obviously, as $ \kappa \to 0+ \Rightarrow f^{**}_h(t) /f^{**}_h(t) \to 1   $ a.e.  Therefore,

$$
\underline{K}(V) \le \lim_{\kappa \to 0+} \left[\frac{|||f^{*}_h|||V}{|||f^{**}_h|||^*V}   \right] =
$$

$$
 \lim_{\kappa \to 0+} \left[\frac{||| 1 - t^{\kappa} |||V}{||| 1 - t^{\kappa}/(\kappa + 1)  |||^*V}   \right] = 1.
$$
 As long as $  \underline{K}(V) \le 1,  $ we   conclude $  \underline{K}(V) = 1.  $ \par

It remains to prove the exactness of right constant $ K(V,H). $\par
 We can suppose without loss of generality the existence  of a positive right continuous decreasing function $  g_0 = g_0(t)  $
from the space $ V $ for which

$$
||g_0||V = 1, \ ||H[g_0]||V = K(V,H)||g_0||V = K(V,H). \eqno(4.5)
$$

 There exists a  function $ f_0 = f_0(x), \ x \in [0,1] $ such that $ f_0^*(t) = g_0(t) $ and following
$$
 ||| f_0 |||_Y =||H[g_0]||V  = K(V,H) \cdot ||g_0||V  =  K(V,H) \cdot ||| f_0 |||^*_Y. \eqno(4.6)
$$
 This completes the proof of theorem 4.1. \par

\vspace{3mm}

 {\it Therefore, we can use for constant $ K(V,H) $ estimate the results of second section.} \par

\vspace{3mm}
 {\bf Example 4.1.} \par
\vspace{3mm}

 Let $ V = L_p(u^p), $ where $ u = u(x) $ is weight function, then $ K(V,H)  $ allows the
following estimate:

$$
K(V,H) \le 2 C_p(u,u), \eqno(4.7)
$$
where $  C_p(u,u) $ is defined in (2.4)-(2.5); and the estimate (4.7) is weakly exact.\par

 As a consequence from theorem 4.1 in the case $ V = L_p(u^p): $

$$
  ||| f |||^*_Y   \le ||| f |||_Y \le 2 C_p(u,u) \cdot ||| f |||^*_Y, \eqno(4.8)
$$

\vspace{3mm}
 {\bf Subexample 4.2.} \par
\vspace{3mm}

   Let in addition $ u(x) = v(x) = x^{\beta}, \ \beta = \const \in [0,1) $ and $ p > 1/(1-\beta); $ then

$$
  ||| f |||^*_Y   \le ||| f |||_Y \le C(\beta) \ \left[\frac{p}{p-1/(1-\beta)} \right] \cdot ||| f |||^*_Y, \eqno(4.9)
$$

\hfill $\Box$ \\
 \bigskip

\section{Generalization on the anisotropic Grand Lebesgue  spaces.}

\vspace{3mm}

{\bf 1. (Ordinary) Grand Lebesgue  spaces.}\par
\vspace{3mm}

  Recently, see  \cite{Fiorenza3}, \cite{Iwaniec1}, \cite{Iwaniec2},
 \cite{Kozachenko1},\cite{Liflyand1}, \cite{Ostrovsky1}, \cite{Ostrovsky2}   etc.
 appear the so-called Grand Lebesgue Spaces $ GLS = G(\psi) =G\psi =
 G(\psi; A,B), \ A,B = \const, A \ge 1, A < B \le \infty, $ spaces consisting
 on all the measurable functions $ f: X \to R $ with finite norms

$$
   ||f||G(\psi) \stackrel{def}{=} \sup_{p \in (A,B)} \left[ |f|_p /\psi(p) \right]. \eqno(5.1)
$$

  Here $ \psi(\cdot) $ is some continuous positive on the {\it open} interval
$ (A,B) $ function such that

$$
     \inf_{p \in (A,B)} \psi(p) > 0, \ \psi(p) = \infty, \ p \notin (A,B).
$$
 We will denote
$$
 \supp (\psi) \stackrel{def}{=} (A,B) = \{p: \psi(p) < \infty, \}
$$

The set of all $ \psi $  functions with support $ \supp (\psi)= (A,B) $ will be
denoted by $ \Psi(A,B). $ \par
  This spaces are rearrangement invariant, see \cite{Bennet1}, and
  are used, for example, in the theory of probability  \cite{Kozachenko1},
  \cite{Ostrovsky1}, \cite{Ostrovsky2}; theory of Partial Differential Equations \cite{Fiorenza2},
  \cite{Iwaniec2};  functional analysis \cite{Fiorenza3}, \cite{Iwaniec1},  \cite{Liflyand1},
  \cite{Ostrovsky2}; theory of Fourier series \cite{Ostrovsky1},
  theory of martingales \cite{Ostrovsky2},mathematical statistics \cite{Sirota2}, \cite{Sirota4}, \cite{Sirota5};
   theory of approximation \cite{Ostrovsky7} etc.\par

 Notice that in the case when $ \psi(\cdot) \in \Psi(A,\infty)  $ and a function
 $ p \to p \cdot \log \psi(p) $ is convex,  then the space
$ G\psi $ coincides with some {\it exponential} Orlicz space. \par
 Conversely, if $ B < \infty, $ then the space $ G\psi(A,B) $ does  not coincides with
 the classical rearrangement invariant spaces: Orlicz, Lorentz, Marcinkiewicz  etc.\par

\vspace{3mm}

{\bf Remark 5.1} If we introduce the {\it discontinuous} function

$$
\psi_r(p) = 1, \ p = r; \psi_r(p) = \infty, \ p \ne r, \ p,r \in (A,B)
$$
and define formally  $ C/\infty = 0, \ C = \const \in R^1, $ then  the norm
in the space $ G(\psi_r) $ coincides with the $ L_r $ norm:

$$
||f||G(\psi_r) = |f|_r.
$$
Thus, the Grand Lebesgue Spaces are direct generalization of the
classical exponential Orlicz's spaces and Lebesgue spaces $ L_r. $ \par

\vspace{3mm}

{\bf Remark 5.2}  The function $ \psi(\cdot) $ may be generated as follows. Let $ \xi = \xi(x)$
be some measurable function: $ \xi: X \to R $ such that $ \exists  (A,B):
1 \le A < B \le \infty, \ \forall p \in (A,B) \ |\xi|_p < \infty. $ Then we can
choose

$$
\psi(p) = \psi_{\xi}(p) = |\xi|_p.
$$

 Analogously let $ \xi(t,\cdot) = \xi(t,x), t \in T, \ T $ is arbitrary set,
be some {\it family } $ F = \{ \xi(t, \cdot) \} $ of the measurable functions:
$ \forall t \in T  \ \xi(t,\cdot): X \to R $ such that
$$
 \exists  (A,B): 1 \le A < B \le \infty, \ \sup_{t \in T} \
|\xi(t, \cdot)|_p < \infty. \eqno(5.2)
$$
 Then we can choose

$$
\psi(p) = \psi_{F}(p) = \sup_{t \in T}|\xi(t,\cdot)|_p.
$$
The function $ \psi_F(p) $ may be called as a {\it natural function} for the family $ F. $
This method was used in the probability theory, more exactly, in
the theory of random fields, see \cite{Ostrovsky1}, chapters 3,4. \par

\vspace{3mm}

{\bf 2. Anisotropic Grand Lebesgue-Riesz spaces.}

\vspace{3mm}
 Let $ Q $ be convex (bounded or not) subset of the set $ \otimes_{j=1}^d [1,\infty]. $
 Let $ \psi = \psi(\vec{p}) $ be continuous in an interior $ Q^0 $ of the set $ Q $
strictly  positive  function such that

$$
\inf_{\vec{p} \in Q^0}  \psi(\vec{p}) > 0; \ \inf_{\vec{p} \notin Q^0}  \psi(\vec{p}) = \infty.
$$

 We denote the set all of such a functions as $ \Psi(Q). $ \par
The  Anisotropic Grand Lebesgue Spaces $ AGLS = AGLS(\psi) $
 space consists on all the measurable functions

$$
f:  \otimes_{j=1}^d X_j \to R
$$
with finite (mixed) norms

$$
||f||AG\psi = \sup_{\vec{p} \in Q^0} \left[ \frac{|f|_{\vec{p}}}{\psi(\vec{p} )} \right]. \eqno(5.3)
$$

 An application  into the theory of multiple Fourier transform of these  spaces see in articles \cite{Benedek1} and
 \cite{Ostrovsky3}, where are considered some problems
of boundedness of singular operators  in (weight) Grand Lebesgue Spaces and  Anisotropic Grand Lebesgue Spaces.
We intend to extend some  results obtained in \cite{Benedek1}, \cite{Ostrovsky3}. \par

\vspace{3mm}

{\bf 3. Hardy's operator in the anisotropic Grand Lebesgue  spaces.}

\vspace{3mm}

 Let $ Q $ be {\it appropriate} for considered problem Hardy's operator estimates domain:
 convex (bounded or not) subset of the set $ \otimes_{j=1}^l [1,\infty]. $
 Let $ \psi = \psi(\vec{p}) $ be continuous in an interior $ Q^0 $ of the set $ Q $
strictly  positive  function such that

$$
\inf_{\vec{p} \in Q^0}  \psi(\vec{p}) > 0; \ \inf_{\vec{p} \notin Q^0}  \psi(\vec{p}) = \infty.
$$

 Let  $ f(\vec{x}) = f(x)  $ be some function  such that the product
$ \vec{x}^{ \vec{\beta} } \ f(\cdot) $  lies in the space $  AG\psi. $
 We denote   the  function $  \vec{q} =\vec{q}(\vec{p})  $ as described before
and denote the inverse function by $  \vec{p} =\vec{p}(\vec{q}).  $ \par

  Define a new function

$$
\nu_R(\vec{q}) = \psi(\vec{p}(\vec{q})) \cdot   K(d; \vec{\alpha}, \vec{\beta}; \vec{p}(\vec{q})),  \eqno(5.4)
$$
 where the parameters $ (\vec{\alpha}, \vec{\beta}, \vec{p}, \vec{q})  $  satisfy the conditions of theorem 3.2.\par
\vspace{3mm}
{\bf Theorem 5.1.}
\vspace{3mm}
$$
|| \vec{x}^{\vec{\alpha}} \ H_d[f] ||AG\nu_R \le
1 \cdot || \vec{x}^{\vec{\beta}} \ f||AG\psi, \eqno(5.5)
$$
{\it where the constant $  "1" $ is the best possible. }\\

\vspace{3mm}

{\bf Proof.}
Let  $ \vec{x}^{ \vec{\beta} } \ f(\cdot) \in AG\psi; $ we can suppose without loss of generality $ ||f||AG\psi = 1. $  This imply that

$$
 || \ \vec{x}^{ \vec{\beta} } \ [f] \ ||_{\vec{p}} \le \psi(\vec{p}).
$$
 We have denoting $ u = \vec{x}^{\vec{\alpha}} \ H_d[f]: $

 $$
 ||u(\cdot)||_q \le   K(d; \vec{\alpha}, \vec{\beta}; \vec{p}(\vec{q}))) \cdot
  \psi( \vec{p}(\vec{q})) \cdot ||f|| AG\psi. \eqno(5.6)
 $$
 As long as the variable $  \vec{p} $ is uniquely defined monotonic function on  $ \vec{q}, $  the inequality (5.5) is
equivalent to the assertion of theorem 5.1.\par
 The {\it exactness} of this estimation is proved in one-dimensional $ d=1 $ in the article
\cite{Ostrovsky3};  the multidimensional case $ d \ge 2 $ provided analogously. \par

\hfill $\Box$ \\
 \bigskip

\section{Concluding remarks.}

\vspace{3mm}

{\bf A.}  The complete investigation of some subclasses of considered here spaces: description of
conjugate (= associate or dual) spaces, conditions of reflexivity  and separability, description of compact subsets,
conditions for absolutely continuity norm, density of simple functions,
boundedness of integral (regular and singular)  operators  with some applications see, e.g. in
\cite{Carro1}, \cite{Carro22}, \cite{Edmund1}, \cite{Edmund2}, \cite{Kaminska1}, \cite{Milman3}, \cite{Opic1}, \cite{Pustylnik1},
\cite{Soria1} etc.\par

\vspace{3mm}

{\bf B.} The fundamental functions $ \Phi_{Y_*}(\delta),  \ \Phi_Y(\delta), \ \delta > 0  $  for the spaces $ (Y, |||\cdot|||Y)  $ and
$ (Y_*, |||\cdot|||Y_*)  $ correspondingly may be calculated as follows. Denote by $ \phi_V(t)  $ the fundamental function of the space
$  (V, ||\cdot||V). $ \par

\vspace{3mm}
{\bf Proposition 6.1.}

$$
\Phi_{Y_*}(\delta) = \phi_V(\delta), \eqno(6.1)
$$

$$
\Phi_{Y}(\delta) = \phi_V(\delta) + || \ I_{(\delta,\mu(X))}(t) \cdot (1/t)  \ ||V.  \eqno(6.2)
$$

\vspace{3mm}
{\bf Proof.} Let $ A $  be measurable subset of $ X $ such that $ \mu(A) = \delta. $ If we denote
$  g(x) = I_A(x), $ then

$$
g^*(t) = I_{[0,\delta]}(t),
$$
hence

$$
\Phi_{Y_*}(\delta)= ||g^*(\cdot)||V = || \ I_{[0,\delta]}(\cdot) \ ||V = \phi_V(\delta).
$$
 Further,

$$
g^{**}(t) = I_{[0,\delta]}(t) + I_{(\delta,\mu(X))}(t) \cdot (1/t),
$$
therefore

$$
\Phi_{Y}(\delta) = || \ I_{[0,\delta]}(\cdot) \ ||V  + || \ I_{(\delta,\mu(X))}(t) \cdot (1/t)  \ ||V   =
\phi_V(\delta) + || \ I_{(\delta,\mu(X))}(t) \cdot (1/t)  \ ||V.
$$

 Recall that the fundamental function play a very important role in the investigation of integral operators
and in the theory of Fourier series and transform, see \cite{Bennet1}, chapter 10.\par

\hfill $\Box$ \\
 \bigskip

{\bf   Acknowledgements.} The authors  would like to thank  prof. S.V.Astashkin, M.M.Milman and L.Maligranda for sending
of Yours remarkable papers and several helpful suggestions. \par

\hfill $\Box$ \\

\end{document}